\documentclass{article}
\usepackage{graphicx}

\input{tcilatex}
\begin{document}

\title{Tilted Cone and Cylinder, Cone and Tilted Sphere}
\author{Mehmet Kirdar}
\maketitle

\begin{abstract}
In this note, we will consider two classical volume problems related to
elliptic integrals. The first problem has a neat formula by means of
elliptic integrals. We remade it with details. In the second problem, we
found a messy formula. On the other hand, it seems to be useful to find a
good approximation for the volume.\vspace{1cm}


\end{abstract}

\section{Introduction}

In this note, I discuss two classical volume problems. The first problem
which, I saw in [2], dates back to 1932, has a neat solution formula by
means of elliptic integrals. I reproduced the formula for the case $k<1$
with some details for elliptic integrals. There is a key identity which also
appeared in the second problem but is not available in [2]. I believe that
Rhodes did this computations somewhere else. His purpose in this article was
Landen transformations but I believe that they are sometimes complications.

The solution for the second problem I found is not very neat. I used
Maclaurin's series expansions of elliptic integrals of the first kind and
the second kind and wrote the solution as an infinite series of
trigonometric integrals. It seems to be useful to find a good approximation
for the volume. I do not know whether this formula was known before. I have
not seen.

I must also mention the beautiful book of Harris Hancock, [1], which helped
me to understand the tricky identities about elliptic integrals.

WolframAlpha helped me a lot during my research. Its abilities are amazing.

\section{Tilted cone and cylinder}

Consider the cylinder\textbf{\ }$x^{2}+y^{2}=1$ and the cone $z=\cot \alpha 
\sqrt{(x-k)^{2}+y^{2}}$, $0\leq k\leq 1.$ We want to find the volume of the
bounded region inside the cylinder, under the cone and above $z=0$. Here $%
\alpha $ is the fixed angle of the cone, the angle between the cone and its
axis, $0\leq \alpha \leq \frac{\pi }{2}.$

Let the origin be $O=(0,0,0),$ the vertex of the cone be $T=(k,0,0)$ and let 
$P=(\cos \theta ,\sin \theta ,0)$, $0\leq \theta \leq 2\pi ,$ be a point one
the unit circle of the $xy$-plane. Let the angle between $TP$ and positive
side of the $x$-axis be\ $\phi $, $0\leq \phi \leq 2\pi $. See [2] for some
figures about this problem. If $TP=R$ then by law of cosines, $R=\sqrt{%
1-k^{2}\sin ^{2}\phi }-k\cos \phi $. The perpendicular from $P$ to $xy$%
-plane cuts the cone with height $R\cot \alpha $. Therefore, the
parameterization of the region in tilted cylindrical coordinates is $0\leq
r\leq R,0\leq \phi \leq 2\pi $ and $0\leq z\leq r\cot \alpha $. And in
tilted coordinates volume differential is $dV=rdrd\phi dz$. With two
successive integrations, the volume integral can be reduced to $V=\dfrac{%
2\cot \alpha }{3}\dint\limits_{0}^{\pi }R^{3}d\phi $.

Now, by putting $R^{3}$ and observing that 
\[
\dint\limits_{0}^{\pi }(-k^{3}\cos ^{3}\phi -3k\cos \phi +3k^{3}\allowbreak
\cos \phi \sin ^{2}\phi )d\phi =\allowbreak 0 
\]
we have%
\[
V=\dfrac{4\cot \alpha }{3}\dint\limits_{0}^{\dfrac{\pi }{2}%
}(3k^{2}+1-4k^{2}\sin ^{2}\phi )\sqrt{1-k^{2}\sin ^{2}\phi }d\phi . 
\]

By the definition of the elliptic integral of the second kind $E(k)$, $V$
becomes 
\[
V=\dfrac{4(3k^{2}+1)\cot \alpha }{3}E(k)-\dfrac{16k^{2}\cot \alpha }{3}%
\dint\limits_{0}^{\dfrac{\pi }{2}}\sin ^{2}\phi \sqrt{1-k^{2}\sin ^{2}\phi }%
d\phi . 
\]

Next, $E_{2}(k)=$ $\dint\limits_{0}^{\dfrac{\pi }{2}}\sin ^{2}\phi \sqrt{%
1-k^{2}\sin ^{2}\phi }d\phi $ must be computed in terms of elliptic
integrals. Let $\Delta =\sqrt{1-k^{2}\sin ^{2}\phi }$, $S=\sin \phi $ and $%
C=\cos \phi $. The tricky identity is

\[
S^{2}\Delta =\frac{2k^{2}-1}{3k^{2}}\Delta +\frac{1-k^{2}}{3k^{2}}\dfrac{1}{%
\Delta }+\left[ -\frac{1}{3}(1-2S^{2})\Delta +\frac{k^{2}}{3}S^{2}C^{2}\frac{%
1}{\Delta }\right] . 
\]

Integrating from $0$ to $\dfrac{\pi }{2}$, since 
\[
\dint\limits_{0}^{\pi /2}\left[ -\frac{1}{3}(1-2S^{2})\Delta +\frac{k^{2}}{3}%
S^{2}C^{2}\frac{1}{\Delta }\right] d\phi =\left[ -\frac{1}{3}SC\Delta \right]
_{0}^{\pi /2}=0, 
\]
we find

\[
E_{2}(k)=\frac{2k^{2}-1}{3k^{2}}E(k)+\frac{1-k^{2}}{3k^{2}}K(k)\text{ }%
(\star ) 
\]%
and

\[
V=\frac{4}{9}\cot \alpha \left[ (k^{2}+7)E(k)+4(k^{2}-1)K(k)\right] 
\]%
where $K(k)=\dint\limits_{0}^{\dfrac{\pi }{2}}\left( 1-k^{2}\sin ^{2}\phi
\right) ^{-\frac{1}{2}}d\phi .$ This formula is obtained in [2]. There,
formula for $k>1$ case is also obtained and then they are combined with a
Landen transformation interpretation.

\section{Cone and tilted sphere}

Let us find the volume of the bounded region between the tilted sphere $%
(x+k)^{2}+y^{2}+z^{2}=1,$ $0\leq k\leq 1$ and the cone $z=\cot \alpha \sqrt{%
x^{2}+y^{2}}.$The set-up of the volume integral is easier than that of the
first problem. So, we can skip figures. The sphere in spherical coordinates
is $\rho ^{2}+2k\rho \cos \theta \sin \phi +k^{2}-1=0$ and the cone is $\phi
=\alpha $. Thus, the volume of the region 
\[
0\leq \rho \leq -k\cos \theta \sin \phi +\sqrt{1-k^{2}+k^{2}\cos ^{2}\theta
\sin ^{2}\phi },\text{ }0\leq \theta \leq 2\pi ,\text{ }0\leq \phi \leq
\alpha 
\]%
is found as

\[
V=\dint\limits_{0}^{2\pi }\dint\limits_{0}^{\alpha }\left\{ 
\begin{array}{c}
(k^{3}\cos \theta \sin ^{2}\phi -k\cos \theta \sin ^{2}\phi -\frac{4}{3}%
k^{3}\allowbreak \cos ^{3}\theta \sin ^{4}\phi )+ \\ 
\left( \frac{4}{3}k^{2}\cos ^{2}\theta \sin ^{3}\phi +\frac{1-k^{2}}{3}\sin
\phi \right) \sqrt{1-k^{2}+k^{2}\cos ^{2}\theta \sin ^{2}\phi }%
\end{array}%
\right\} d\phi d\theta . 
\]

Since%
\[
\dint\limits_{0}^{2\pi }(k^{3}\cos \theta \sin ^{2}\phi -k\cos \theta \sin
^{2}\phi -\frac{4}{3}k^{3}\allowbreak \cos ^{3}\theta \sin ^{4}\phi )d\theta
=\allowbreak 0 
\]
and due to symmetry, we find

\[
V=\dint\limits_{0}^{\alpha }\dint\limits_{0}^{\dfrac{\pi }{2}}\left( \frac{16%
}{3}k^{2}\cos ^{2}\theta \sin ^{3}\phi +\frac{4}{3}(1-k^{2})\sin \phi
\right) \sqrt{1-k^{2}+k^{2}\cos ^{2}\theta \sin ^{2}\phi }d\theta d\phi . 
\]

Let us define 
\[
K=\dfrac{k\sin \phi }{\sqrt{1-k^{2}\cos ^{2}\phi }} 
\]
and thus,

\begin{eqnarray*}
V &=&\dint\limits_{0}^{\alpha }\left( \frac{16}{3}k^{2}\sin ^{3}\phi +\frac{4%
}{3}(1-k^{2})\sin \phi \right) \sqrt{1-k^{2}\cos ^{2}\phi }E(K)d\phi \\
&&-\dint\limits_{0}^{\alpha }\text{ }\frac{16}{3}k^{2}\sin ^{3}\phi \sqrt{%
1-k^{2}\cos ^{2}\phi }\left( \dint\limits_{0}^{\dfrac{\pi }{2}}\sin
^{2}\theta \sqrt{1-K^{2}\sin ^{2}\theta }d\theta \right) d\phi .
\end{eqnarray*}

By using the star identity, $(\star )$ of the first problem, it can be
written as

\[
V=\frac{4}{9}\dint\limits_{0}^{\alpha }(8k^{2}\sin ^{3}\phi +7(1-k^{2})\sin
\phi )\sqrt{1-k^{2}\cos ^{2}\phi }E(K)d\phi -\frac{16}{9}\dint\limits_{0}^{%
\alpha }\sin \phi \sqrt{1-k^{2}\cos ^{2}\phi }K(K)d\phi . 
\]

We can now insert infinite series of $E(K)$ and $K(K)$ and do term by term
integration to obtain a formula which involves trigonometric integrals.

Let us recall that $E(K)=\dfrac{\pi }{2}\dsum\limits_{n=0}\dfrac{c_{n}}{1-2n}%
K^{2n}$ and $K(K)=\dfrac{\pi }{2}\dsum\limits_{n=0}c_{n}K^{2n}$ where $%
c_{n}=\left( \dfrac{(2n)!}{2^{2n}(n!)^{2}}\right) ^{2}.$ Putting these and $%
K $ in the last equation we obtain

\[
V=\frac{2\pi }{9}\dsum\limits_{n=0}^{\infty }\frac{c_{n}k^{2n}}{1-2n}%
\dint\limits_{0}^{\alpha }\frac{8k^{2}\sin ^{2n+3}\phi +(3-7k^{2}+8n)\sin
^{2n+1}\phi }{(1-k^{2}\cos ^{2}\phi )^{n-\frac{1}{2}}}d\phi . 
\]

The zeroth term of the series gives the following approximation of the
volume for small $k$: 
\[
\frac{2\pi }{9}\left( (1+2k^{2})\sqrt{1-k^{2}}-\cos \alpha
(1+4k^{2}-2k^{2}\cos \alpha )\sqrt{1-k^{2}\cos ^{2}\alpha }+(2-3k^{2})\frac{%
\arcsin k-\arcsin (k\cos \alpha )}{k}\right) . 
\]%
This gives the exact value $\dfrac{2\pi }{3}(1-\cos \alpha )$ in the limit $%
k\rightarrow 0.$

\section{Acknowledgement}

I thank Paul Bracken who introduced me with AGM\ and elliptic integrals.

$\allowbreak $

\end{document}